\documentclass[11pt,leqno]{article}

\usepackage{amsmath,amsfonts,amscd,amssymb,theorem}

\long\def\comment#1\endcomment{}

\comment
\endcomment

\makeatletter
\begingroup
\gdef\th@dotted{\normalfont\itshape
  \def\@begintheorem##1##2{%
        \item[\hskip\labelsep \theorem@headerfont ##1\ ##2.]}%
\def\@opargbegintheorem##1##2##3{%
   \item[\hskip\labelsep \theorem@headerfont ##1\ ##2\ (##3).]}}
\endgroup
\makeatother

\theoremstyle{dotted}

\newtheorem{theorem}{Theorem}[section]
\newtheorem{lemma}[theorem]{Lemma}

\newtheorem{prop}[theorem]{Proposition}
\newtheorem{corr}[theorem]{Corollary}

\makeatletter
\begingroup
\gdef\th@upshape{\normalfont
  \def\@begintheorem##1##2{%
        \item[\hskip\labelsep \theorem@headerfont ##1\ ##2.]}%
\def\@opargbegintheorem##1##2##3{%
   \item[\hskip\labelsep \theorem@headerfont ##1\ ##2\ (##3).]}}
\endgroup
\makeatother

\theoremstyle{upshape}

\newtheorem{defn}[theorem]{Definition}
\newtheorem{remark}[theorem]{Remark}

\makeatletter
\renewcommand{\subsection}{\@startsection{subsection}{2}{0pt}{-3ex
plus -1ex minus -0.2ex}{-2mm plus -0pt minus
-2pt}{\normalfont\bfseries}} 
\renewcommand{\subsubsection}{\@startsection{subsubsection}{3}{0pt}{-3ex
plus -1ex minus -0.2ex}{-2mm plus -0pt minus
-2pt}{\normalfont\bfseries}} 
\makeatother

\makeatletter
\@addtoreset{equation}{section}
\makeatother
\renewcommand{\theequation}{\thesection.\arabic{equation}}

\newcommand{\cntrct}                % contraction with a vector field
{\hspace{2pt}\raisebox{1pt}{\text{$\lrcorner$}}\hspace{2pt}}

\newcommand{\proof}[1][Proof.]{\smallskip\noindent{\em #1}}
\def\endproof{\hfill\ensuremath{\square}\par\medskip}

\def\eqref#1{\thetag{\ref{#1}}}

\let\latexref=\ref
\def\ref#1{{\normalfont{\latexref{#1}}}}

\newcommand{\wt}{\widetilde}

\newcommand{\idot}{{\:\raisebox{1pt}{\text{\circle*{1.5}}}}}
\newcommand{\eps}{\varepsilon}
\renewcommand{\phi}{\varphi}

\newcommand{\Hom}{\operatorname{Hom}}
\newcommand{\Ext}{\operatorname{Ext}}

\renewcommand{\Im}{\operatorname{Im}}

\newcommand{\id}{\operatorname{\sf id}}

\newcommand{\gr}{\operatorname{\sf gr}}

\newcommand{\D}{{\cal D}}
\newcommand{\C}{{\cal C}}

\newcommand{\hash}{\sharp}

\newcommand{\End}{{\operatorname{End}}}

\newcommand{\ppt}{{\sf pt}}

\newcommand{\cchar}{\operatorname{\sf char}}

\newcommand{\Z}{{\mathbb Z}}

\newcommand{\hocolim}{\operatorname{\sf hocolim}}

\newcommand{\Ss}{\mathbb{S}}
\newcommand{\Nn}{\mathbb{N}}

\newcommand{\Dalg}{\operatorname{\mathcal{D}\mathcal{A}\mathit{l}\mathit{g}}}
\newcommand{\Dcomm}{\operatorname{\mathcal{D}\mathcal{C}\mathit{o}
    \mathit{m}\mathit{m}}}

\newcommand{\eCat}{\operatorname{\mathcal{C}\mathit{a}\mathit{t}}}

\newcommand{\Q}{\mathbb{Q}}

\newcommand{\Der}{\operatorname{Der}}
\newcommand{\Iso}{\operatorname{Iso}}

\newcommand{\Ind}{\operatorname{Ind}}

\newcommand{\F}{\mathbb{F}}

\newcommand{\bR}{R_0}

\title{Spectral algebras and non-commutative Hodge-to-de Rham degeneration}

\author{D. Kaledin, A. Konovalov and K. Magidson\thanks{All
    authors were partially supported by Basis Foundation grant
    18-1-6-95-1, Leader (Math). D.K. and A.K. were also partially supported
    by the HSE University Basic Research Program, Russian Academic Excellence
    Project '5-100'}}

\date{{\em To the blessed memory of I.R. Shafarevich}}

\begin{document}

\maketitle

\tableofcontents

\section*{Introduction.}

For any DG algebra $A$ over a field $K$, one has the {\em
  Hochschild-to-cyclic}, or {\em Hodge-to-de Rham} spectral sequence
$$
HH_\idot(A)((u)) \Rightarrow HP_\idot(A), \qquad \deg u = 2
$$
relating Hochschild and periodic cyclic homology of the algebra
$A$. It has been conjectured by Kontsevich and Soibelman \cite{KS}
that if $\cchar K = 0$ and $A$ is smooth and proper, the spectral
sequence degenerates. The conjecture has been proved under some
restrictions in \cite{Ka1}, and in full generality in
\cite{Ka2}. Recently, a slightly different proof was given by
A. Mathew in \cite{M}.

This paper arose as an attempt to generalize these results to other
settings of interest for applications (for example, to
$\Z/2\Z$-graded DG algebras). As of now, we did not succeed;
however, we think that we can at least streamline and clarify the
original proofs of \cite{Ka1}, \cite{Ka2}. This is the subject of
the present paper.

While the degeneration statement itself is purely homological, all
the proofs use stable homotopy theory. This is quite explicit in
\cite{Ka1}, even more explicit in Mathew's proof, and implicitly
present also in \cite{Ka2} (actually, it was deliberately hidden so
as to accomodate the readers who do not like topology). The main
reason why topology could possibly help can be summarized as
follows:
\begin{itemize}
\item If an algebra $A$ is smooth and proper over $K$, then its
  Hodge-to-de Rham spectral sequence consists of finite-dimensional
  $K$-vector spaces, so by the standard criterion of Deligne, it
  degenerates if and only if the first page is abstractly isomorphic
  to the last one. More generally, Hochschild Homology $HH(A/R)$
  exists for an algebra $A$ over any commutative ring spectrum $R$,
  and we can ask whether there exists an isomorphism
\begin{equation}\def\theequation{*}
HH(A/R) \otimes_R R^{tS^1} \cong HP(A/R),
\end{equation}
where $R^{tS^1}$ stands for the Tate homotopy fixed points of the
spectrum $R$ with respect to the trivial action of the circle $S^1$,
and $HP(A/R) = HH(A/R)^{tS^1}$ are the Tate fixed points of
$HH(A/R)$ with respect to the standard circle action. The homotopy
groups $\pi_\idot(R^{tS^1})$ can be computed by the
Atiyah-Hirzeburch spectral sequence that starts at
$\pi_\idot(R)((u))$. If $R$ is orientable -- for example, if $R$ is
a usual commutative ring -- then the sequence degenerates, so that
$R^{tS^1} \cong R((u))$. But in general, it does not have to, so
that $R^{tS^1}$ can be smaller that $R((u))$. Under favourable
circumstances, it can become so small that \thetag{*} exists for
trivial reasons.
\end{itemize}
In practice, we do not know whether these ``favourable
circumstances'' really occur. However, if one considers a cyclic
subgroup $C_p \subset S^1$ of some prime order $p$, then a striking
result known as the Segal Conjecture shows that for the sphere
spectrum $\Ss$, the Tate fixed point spectrum $\Ss^{tC_p}$ is simply
the $p$-completion $\Ss_p$ -- that is, it is as small as it could
possibly be (in particular, it is connective). This suggests that
one should consider separately all primes, and prove the theorem by
reducing the statement at each prime $p$ to a statement about the
Tate fixed points $HH(A/R)^{tC_p}$ that would follow from the Segal
Conjecture.

If one cuts down to the point, then this is exactly what happens in
\cite{Ka1} and \cite{Ka2}. Formally, the argument replicates the
classic proof of the commutative Hodge-to-de Rham degeneration of
Deligne and Illusie \cite{DI}, and it works by reduction to positive
characteristic. The reduction is achieved by a beautiful theorem of
B. T\"oen stating that $A \cong A_R \otimes_R K$ for some smooth and
proper DG algebra $A_R$ finitely generated subring $R \subset K$
smooth over $\Z$. Then for each residue field $k$ of some positive
characteristic $p$, one needs to prove degeneration for $A_k = A_R
\otimes_R k$. While in general, Hodge-to-de Rham degeneration in
positive characteristic is false, it still holds under additional
assumptions. In \cite{Ka1}, \cite{Ka2}, the assumptions are that $A$
lifts to the second Witt vectors ring $W_2(k)$, and that Hochschild
cohomology $HH^i(A)$ vanishes for $i \geq 2p$. What the second
assumption really means though, explicitly in \cite{Ka1} and
implicitly in \cite{Ka2}, is that $A$ can be lifted to an algebra
over a certain ring spectrum, a topological counterpart of the ring
$W_2(k)$. Degeneration is then due to some very truncated version of
the Segal Conjecture for the group $C_p$ proved essentially by hand.

Mathew in \cite{M} has similar assumptions but with Hochschild
cohomology replaced by Hochschild homology, and this is because his
strategy is different: instead of lifting a $k$-algebra $A$ to a
spectrum, he considers it as a spectrum as it is, and then uses deep
results about Topological Hochschild Homology for his proof. It is
hard to see how this can be improved, but in retrospect, it is
obvious what can be done with \cite{Ka2}. Instead of first
restricting our algebra $A$ to a ring $R \subset K$, then localizing
$R$ to insure that all its residue fields $k$ are good enough, and
then lifting each reduction $A_k$ to an algebra over a ring spectrum
by obstruction theory, one should directly restrict $A$ to an
approriate ring spectrum $R$, so that there is no need to lift and
no conditions to impose. This is the argument that we sketch in this
paper.

\medskip

One obvious problem with this streamlined argument is that it really
has to be done topologically, and one needs an appropriate
technology for that. It is more-or-less clear by now that ideally,
one would like to have some model-independent formalism of
``enhanced categories'', both stable and unstable, and this
formalism should be equipped with a concise and convenient toolkit
sufficient for practical applications. At present, the only existing
formalism is that of $\infty$-categories in the sense of J. Lurie,
and that is not model-independent (instead of choosing a category of
models, you have to choose a model of your category). What is worse,
it does not differentiate cleanly between the model-dependent and
model-independent parts, and cannot be used as a black box. There is
no convenient toolkit --- on the contrary, a rigorous paper written
in the $\infty$-categorical language has to rely on several thousand
pages of Lurie's foundational work, and to give precise references
at every second line. In principle, it is possible to do this; a
perfect example is the recent paper \cite{NS}. However, it seems
that the widespread practice these days is to not to do this, and
rely instead on the reader's conjectural capability to fill in all
the missing details.

We emphasize that this is very bad practice that is certain to lead
to disaster, and we choose to follow suit. Our justification is that
after all, the Degeneration Theorem has been already proved. Our
goal is to explain the proof and show how it can be improved, not to
re-do it with complete rigor. Conversely, having a concrete,
detailed and non-trivial application can show what needs to be a
part of any usable future toolkit, and possibly help develop it. To
emphasize the provisional nature of our results, we speak of
enhanced categories and functors instead of $\infty$-categories, and
we state clearly that what we have in the paper is no more than a
sketch.

\subsection*{Acknowledgement.} We are grateful to A. Efimov,
A. Fonarev, L. Hesselholt, Th. Nikolaus and A. Prihodko for useful
discussions, and to MSRI where part of this work was done. We are
especially grateful to A. Mathew for generously sharing his insights
and expertise, and in particular, for helping us with the (sketch of
the) proof of Proposition~\ref{sm.prop}.

\section{Preliminaries.}

\subsection{Enhanced categories.}

For any enhanced category $\C$, we denote by $\pi_0(\C)$ its
truncation to an ordinary category. An enhanced functor $\gamma:\C
\to \C'$ induces a functor $\pi_0(\gamma):\pi_0(\C) \to \pi_0(\C')$
that we will denote simply by $\gamma$ if there is no danger of
confusion. For any enhanced category $\C$ and small category $I$,
enhanced functors from $I$ to $\C$ form an enhanced category
$\C^I$. We have a natural conservative comparison functor
\begin{equation}\label{pi.eq}
\pi_0:\pi_0(\C^I) \to \pi_0(\C)^I,
\end{equation}
and if $I = \Nn$ is the totally ordered set of positive integers
considered as a small category in the usual way, then \eqref{pi.eq}
is essentially surjective and full. A functor $\gamma:I_0 \to I_1$
induces an enhanced pullback functor $\gamma^*:\C^{I_1} \to
\C^{I_0}$. An enhanced category $\C$ is {\em cocomplete} if for any
small $I$, the pullback functor $\tau^*:\C \to \C^I$ induced by the
projection $\tau:I \to \ppt$ to the point category $\ppt$ admits a
left-adjoint enhanced functor $\hocolim_I:\C^I \to \C$. An object $c
\in \C$ in a cocomplete enhanced category $\C$ is {\em compact} if
the Yoneda enhanced functor $\Hom(c,-)$ commutes with $\hocolim_I$
for any small filtered $I$. A cocomplete enhanced category $\C$ is
{\em compactly generated} if the full enhanced subcategory $\C^{pf}
\subset \C$ spanned by compact objects is small, and for any object
$c \in \C$, we have $c \cong \hocolim_Ic_\idot$ for an enhanced
functor $c_\idot:I \to \C^{pf}$ from a filtered small category
$I$. Any small enhanced category $\C$ canonically embeds as a fully
faithful enhanced subcategory into its {\em $\Ind$-completion}
$\Ind(\C)$; this is a cocomplete compactly generated enhanced
category, and any $c \in \C$ is compact in $\Ind(\C) \supset \C$.

Small enhanced categories themselves form an enhanced category
$\eCat$. This category is cocomplete. The full enhanced subcategory
$\eCat^{\leq 1} \subset \eCat$ spanned by ordinary small categories
is closed under filtered homotopy colimits (but not under all
colimits), and truncation defines an enhanced functor $\pi_0:\eCat
\to \eCat^{\leq 1}$ left-adjoint to the embedding. The functor
$\pi_0$ commutes with filtered homotopy colimits, and filtered
homotopy colimits in $\eCat^{\leq 1}$ are the classical $2$-colimits
of ordinary categories.

We will say that an enhanced category $\C$ is {\em Karoubi-closed}
if so is its truncation $\pi_0(\C)$. The following useful lemma is
essentially due to B. T\"oen.

\begin{lemma}\label{toen.le}
Assume given an enhanced functor $\gamma:\C' \to \C$ between
cocomplete enhanced categories that preserves filtered homotopy
colimits, and assume that $\pi_0(\gamma)$ is conservative and $\C$
is Karoubi-closed. Then $\C'$ is Karoubi-closed.
\end{lemma}

\proof{} Assume given an object $c \in \C'$ and a idempotent
endomorphism $p:c \to c$ in $\pi_0(\C')$, $p^2=p$. Let $C:\Nn \to
\C'$ be the constant enhanced functor with value $c$, and consider
the functor $C(p)_0:\Nn \to \pi_0(\C')$ sending any integer $n \in
\Nn$ to $c$, with transition maps $C(p)_0(n) \to C(p)_0(n+1)$ equal
to $p$. Let $B_0:C(p)_0 \to \pi_0(C)$ be the map equal to $p$ at any
$n \in \Nn$. Since the functor \eqref{pi.eq} is essentially
surjective and full for $I = \Nn$, we can lift $C(p)_0$ to an
enhanced functor $C(p):\Nn \to \C$, $\pi_0(C(p)) \cong C(p)_0$, and
$B_0$ lifts to a map $B:C(p) \to C(p)$ of enhanced functors. By
adjunction, the isomorphism $c \cong C(p)(0)$ induces a map $A:C \to
C(p)$. Since $\C'$ is cocomplete, $\hocolim_{\Nn}$ exists and is
functorial, and if we let $c(p) = \hocolim_{\Nn}C(p)$, then $A$ and
$B$ induce maps
$$
a:c = \hocolim_{\Nn}C \to c(p), \qquad b:c(p) \to c.
$$
Again by adjunction, we have $b \circ a = p$. If the idempotent $p$
does have an image $c'$ --- that is, we have $c' \in \C'$ and maps
$a':c \to c'$, $b':c' \to c$ such that $b' \circ a' = p$ and $a'
\circ b' = \id$ in $\pi_0(\C')$ --- then one easily checks that the
composition $a \circ b':c' \to c(p)$ is an isomorphism, so that $a
\circ b = \id$ by the uniqueness of idempotent images. If not, then
since $\gamma$ commutes with filtered homotopy colimits and $\C$ is
Karoubi-closed, we at least see that $\gamma(a \circ b) = \gamma(a)
\circ \gamma(b)=\id$ in $\pi_0(\C)$. But since $\gamma$ is
conservative, this implies that $a \circ b$ is invertible, and then
$$
(a \circ b)^3 = a \circ (b \circ a)^2 \circ b = a \circ p^2 \circ b
= a \circ p \circ b = (a \circ b)^2,
$$
so that $a \circ b = \id$.
\endproof

\subsection{Spectral algebras.}

We denote by $\D(\Ss)$ the stable enhanced category of spectra. It
is cocomplete, compactly generated and Karoubi-closed (the latter is
slightly non-trivial since e.g.\ the enhanced category of unpointed
homotopy types is not). It also carries a natural structure of a
symmetric monoidal enhanced category, and the enhanced categories
$\Dalg(\Ss)$, $\Dcomm(\Ss)$ of $E_1$ resp.\ $E_\infty$-algebras in
$\D(\Ss)$ are also cocomplete. The stable enhanced category
$\D(\Ss)$ --- or strictly speaking, its triangulated trunctation
$\pi_0(\D(\Ss))$ --- carries a natural $t$-structure, with
$\D^{\leq 0}(\Ss) \subset \D(\Ss)$ consisting of connective
spectra, and a spectrum is {\em discrete} if it lies in the heart of
this natural $t$-structure. Sending $E$ to $\pi_0(E)$ identifies the
heart with the category of abelian groups. An $E_\infty$-algebra in
$\D(\Ss)$ that is discrete is the same thing as unital associative
commutative ring.

For any positive integer $N$, we let $\Ss(N^{-1})$ be the
localization of the sphere $\Ss$ in $N$. We note that we have
\begin{equation}\label{Q.co}
\Q \cong \hocolim_N\Ss(N^{-1}),
\end{equation}
where the colimit is taken with respect to the divisibility order,
and $\Q$ is the field of rationals considered as a discrete
$E_\infty$-algebra in $\D(\Ss)$.

For any $E_1$-algebra $A \in \Dalg(\Ss)$, we have the cocomplete
stable enhanced category $\D(A)$ of left $A$-modules, and for any
$E_\infty$-algebra $R$ in $\Dcomm(\Ss)$, we have the cocomplete
stable symmetric monoidal enhanced category $\D(R)$ of
$R$-modu\-les, and the cocomplete symmetric monoidal enhanced
categories $\Dalg(R)$, $\Dcomm(R)$ of $E_1$
resp.\ $E_\infty$-algebras in $\D(R)$. A map $R \to R'$ between
$E_1$ or $E_\infty$-algebras induces an adjoint pair of the tensor
product functor $\D(R) \to \D(R')$, $M \mapsto R' \otimes_R M$, and
the restriction functor $\D(R') \to \D(R)$. In the $E_\infty$-case,
the tensor product functor is symmetric monoidal, while the
restriction functor is lax symmetric monoidal by adjunction (in the
$\infty$-categorical setup, this is \cite[Corollary 7.3.2.7]{lu17});
therefore they induce adjoint pairs of functors between $\Dalg(R)$
and $\Dalg(R')$, and between $\Dcomm(R)$ and $\Dcomm(R')$. In all
these adjoint pairs, the restriction functor commutes with filtered
colimits, so that by adjunction, the tensor product functor sends
compact objects to compact objects.

The enhanced category $\D(R)$ is compactly generated but there is
more. Namely, the forgetful functor $\D(R) \to \D(\Ss)$ has a
left-adjoint free module functor $F:\D(\Ss) \to \D(R)$, $F(V) = V
\otimes_{\Ss} R$, and an object $M \in \D(R)$ is {\em finitely
  presented} if it is a finite homotopy colimit of objects of the
form $F(E)$, $E \in \D(\Ss)^{pf}$. Then any object in $\D(R)$ is a
filtered homotopy colimit of finitely presented objects. Since
filtered colimits commute with finite limits, any finitely presented
object in $\D(R)$ is compact, so are its retracts, and conversely,
since an isomorphism $M \cong \hocolim_IM_i$ with filtered $I$ and
compact $M$ must factor through some $M_i$, a compact object is a
retract of a finitely presented one. Exactly the same holds for
$\Dalg(R)$, $\Dcomm(R)$, and $\D(A)$ for any $A \in
\Dalg(R)$. Moreover, the forgetful functor is conservative and
commutes with filtered homotopy colimits, so that $\D(R)$ is
Karoubi-closed by Lemma~\ref{toen.le}, and again, the same holds for
$\Dalg(R)$, $\Dcomm(R)$ and $\D(A)$. Furthermore, we have the full
subcategories $\D^{\leq 0}(\Ss) \subset \D(\Ss)$, $\Dalg^{\leq
  0}(\Ss) \subset \Dalg(\Ss)$, $\Dcomm^{\leq 0}(\Ss) \subset
\Dcomm(\Ss)$ spanned by connective spectra, these are also compactly
generated, and so are $\D^{\leq 0}(R) \subset \D(R)$, $\Dalg^{\leq
  0}(R) \subset \Dalg(R)$, $\Dcomm^{\leq 0}(R) \subset \Dcomm(R)$
for any connective $E_\infty$-algebra $R \in \Dcomm^{\leq 0}(\Ss)$.

\begin{remark}
Compact objects in $\D(R)$ are also known as {\em perfect
  $R$-modu\-les}; this explains our notation (although one usually
writes $\D^{pf}(R)$ instead of $\D(R)^{pf}$). For algebras, there is
no standard terminology. T\"oen calls compact algebras
``homotopically finitely presented''.
\end{remark}

For any $R \in \Dcomm(\Ss)$ and $A \in \Dalg(R)$, the cocomplete
enhanced category $\D(A)$ coincides with the $\Ind$-completion
$\Ind(\D(A)^{pf})$ of its full subcategory of compact objects. Aside
from compactness, there is another useful finiteness conditions one
can impose on $A$-modules: an $A$-module $M \in \D(A)$ is {\em
  coherent} if it is compact as an object in $\D(R)$. We note that
unlike compactness, the property of being coherent is preserved by
restriction via an algebra map. In fact, any compact object in
$\D(R)$ is dualizable, so that we have the endomorphism algebra
$\End_R(M) \in \Dalg(R)$, and $M$ is canonically an
$\End_R(M)$-module. Then $M$ is tautologically coherent over
$\End_R(M)$, and any structure of an $A$-module on $M$ is induced
from this canonical $\End_R(M)$-module structure by restriction via
an action map $a:A \to \End_R(M)$ in $\Dalg(R)$.  We denote by
$\D(A)^{coh} \subset \D(A)$ the full subcategory spanned by coherent
modules, and we note that its $\Ind$-completion $\Ind(\D(A)^{coh})$
is in general different from $\D(R)$.

For any $E_\infty$-algebra $R \in \Dcomm(\Ss)$, sending an $E_1$ or
an $E_\infty$-algebra $R'$ over $R$ to the enhanced category
$\D(R')^{pf}$ of compact $R$-modules gives enhanced functors
\begin{equation}\label{mod.f.eq}
\D^{pf}:\Dalg(R),\Dcomm(R) \to \eCat,
\end{equation}
while sending $R \in \Dcomm(\Ss)$ to $\Dalg(R)^{pf}$,
$\Dcomm(R)^{pf}$ give functors
\begin{equation}\label{alg.f.eq}
\Dalg^{pf},\Dcomm^{pf}:\Dcomm(\Ss) \to \eCat.
\end{equation}
We will need the following fundamental fact.

\begin{prop}\label{filt.prop}
The enhanced functors \eqref{mod.f.eq} and \eqref{alg.f.eq} commute
with filtered homotopy colimits.
\end{prop}

\proof[Outline of a proof.] The argument is the same in all cases.
For finitely presented objects $M = \hocolim_IF(E_i)$, $I$ finite,
the proof is a straightforward induction on the cardinality of
$I$. In general, use the characterization of compact objects as
retracts of finitely presented ones, and observe that as we have
already proved, the necessary retractions must also appear at some
finite level.
\endproof

\section{Formal smoothness.}

For any $E_\infty$-algebra $A \in \Dcomm(\Ss)$, any
$E_\infty$-algebra $R \in \Dcomm(A)$, and any $R$-module $M \in
\D(R)$, we have the split square-zero extension $R \oplus M \in
\D(A)$ of $R$ by $M$, and derivations from $R$ to $M$ are splittings
$R \to R \oplus M$ of the augmentation map $R \oplus M \to
R$. Derivations form an enhanced functor $\Der:\D(R) \to \D(\Ss)$
that is representable by the {\em cotangent module} $\Omega(R/A) \in
\D(R)$. If $R$ is compact in $\Dcomm(A)$, then $\Omega(R/A)$ is
compact in $\D(R)$. The same module also controls non-split
square-zero extensions. In particular, if there are no maps from
$\Omega(R/A)$ to the homological shift $M[1]$ of some $M \in \D(R)$,
then any square-zero extension
$$
\begin{CD}
M @>>> R' @>>> R
\end{CD}
$$
in $\Dcomm(A)$ admits a splitting $R \to R'$. The cotangent module
$\Omega(-/A)$ is functorial in the appropriate sense, and commutes
with filtered colimits: for any enhanced functor $R_\idot:I \to
\Dcomm(A)$ with small filtered $I$ and $R = \hocolim_IR_\idot$, we
have an enhanced functor from $I$ to $\D(R)$ with values
$\Omega(R_i/A) \otimes_{R_i} R$, and a natural isomorphism
\begin{equation}\label{om.hoco}
\Omega(R/A) \cong \hocolim_I\Omega(R_\idot/A) \otimes_{R_\idot} R.
\end{equation}

\begin{remark}
In the $\infty$-categorical setting, the sketch above corresponds to
\cite[Section 7.3,7.4]{lu17}, but for some reason, the logic there
is reversed: instead of first defining square-zero extensions,
e.g.\ by considering the natural symmetric monoidal structure on the
filtered version of $\D(\Ss)$, Lurie first defines derivations. The
end result is the same.
\end{remark}

For any $E_\infty$-algebra $R \in \Dcomm(\Ss)$ and any set $S$, we
have the free $R$-module $R[S] \in \D(R)$ generated by $S$. We say
that $M \in \D(R)$ is {\em projective} if it is a retract of a free
$R$-module $R[S]$, and {\em finitely generated projective} if $S$
can be chosen to be finite. A finitely generated projective module
is compact, and conversely, a compact projective module is finitely
generated.

\begin{defn}\label{sm.def}
For any connective $A \in \Dcomm(\Ss)$, an $E_\infty$-algebra $R \in
\Dcomm(A)$ is {\em formally smooth} if it is connective, compact in
$\Dcomm(A)$, and $\Omega(R/A)$ is a projective $R$-module.
\end{defn}

If $A=\Q$ is the field of rationals, then $\D(\Q)$ is the derived
category of complexes of $\Q$-vector spaces, $\Dcomm(\Q)$ is the
category of commutative DG algebras over $\Q$, and $A \in
\Dcomm(\Q)$ is formally smooth iff it is a finitely generated smooth
$\Q$-algebra placed in the homological degree $0$. Over $\Ss$,
formally smooth algebras are not that easy to describe. However,
observe that if $R \in \Dcomm(\Ss)^{pf}$ is formally smooth, then
$\pi_0(R)$ is at least a finitely generated commutative ring.

\begin{prop}\label{sm.prop}
For any field $K$ of characteristic $0$, there exists an isomorphism
$K \cong \hocolim_IR_\idot$ for some small filtered $I$ and an
enhanced functor $R_\idot:I \to \Dcomm(\Ss)^{pf}$ whose values
$R_i$, $i \in I$ are formally smooth in the sense of
Definition~\ref{sm.def}.
\end{prop}

\proof{} Since $\Dcomm(\Ss)$ is compactly generated, we may assume
that $K \cong \hocolim_IR_\idot$ for some small filtered $I$ and
$R_\idot:I \to \Dcomm(\Ss)^{pf}$. Moreover, since $K$ is connective
and $\Dcomm^{\leq 0}(\Ss)$ is also compactly generated, we may
assume that all the $R_i$ are connective. What we need to check is
that one can arrange for them to be formally smooth. For this, it
suffices to show that any map $r:R \to K$ from a compact connective
$R \in \D(\Ss)^{pf}$ factors through a formally smooth
$E_\infty$-algebra $C$.

Indeed, any finitely generated subring $C_0 \subset K$ lies in a
finitely generated smooth $\Q$-algebra $C \subset K$. Since $R$ is
connective, we have the augmentation map $a:R \to \pi_0(R)$, and $r
= b \circ a$ for some map $b:\pi_0(R) \to K$. Then $\pi_0(R)$ is
finitely generated, and taking $C_0 = \Im b$, we see that $r$
factors through a finitely generated smooth $\Q$-subalgebra $C
\subset K$. But then, by Proposition~\ref{filt.prop}, $\Dcomm^{pf}$
commutes with filtered homotopy colimits, and in particular, it
commutes with the colimit \eqref{Q.co}. Thus $C = C_N
\otimes_{\Ss[N^{-1}]} \Q$ for some positive integer $N$. Moreover,
since $C$ is formally smooth over $\Q$, we have maps $a:\Omega(C/\Q)
\to \Q[S]$, $b:\Q[S] \to \Omega(C/\Q)$, $b \circ a = \id$ for some
finite set $S$, and again by Proposition~\ref{filt.prop}, we can
assume after enlarging $N$ that both are induced by maps
$a_N:\Omega(C_N/\Ss(N^{-1})) \to \Ss(N^{-1})[S]$,
$b_N:\Ss(N^{-1})[S] \to \Omega(C_N/\Ss(N^{-1}))$ such that $b_N
\circ a_N = \id$. Therefore $C_N$ is formally smooth over
$\Ss(N^{-1})$, and then also over $\Ss$ since $\Ss \to \Ss(N^{-1})$
is a localization. Finally, since $R$ is compact, we can again
enlarge $N$ so that the map $R \to C \cong \hocolim_NC_N$ factors
through $C_N$, and this finishes the proof.
\endproof

Now, for any prime $p$, denote by $\Ss_p \in \Dcomm(\Ss)$ the
$p$-completion of the sphere $\Ss$, with its natural map $\Ss_p \to
\F_p$, and for any power $q=p^n$ of $p$, let $\Ss_q$ be the $n$-fold
Galois extension of $\Ss_p$, with its map $\Ss_q \to \F_q$ (since
$\F_q$ is \'etale over $\F_p$, the cotangent complex
$\Omega(\F_q/\F_p)$ vanishes, so that $\Ss_q$ exists and is unique).

\begin{lemma}\label{fm.le}
Assume given an algebra $R \in \D(\Ss)$ formally smooth in the sense
of Definition~\ref{sm.def}. Then for any finite field $k = \F_q$,
any map $a:R \to k$ factors through the canonical map $\Ss_q \to k$.
\end{lemma}

\proof{} The completed sphere $S=\Ss_q$ is the homotopy limit of an
enhaced functor $S_\idot:\Nn^o \to \Dcomm(\Ss)$, $n \geq 1$, where
$S_1 = k$, and each $S_{n+1}$ is a square-zero extension of $S_n$ by
a connective $k$-module $M \in \D(k)$. Since \eqref{pi.eq} is full
and essentially surjective for $I=\Nn$, it suffices to extend
$a_1=a:R \to S_1=k$ to a compatible system of factorizations $a_n:R
\to S_n$, $n \geq 2$. This can be done by induction: at each step,
the obstruction to lifting $a_n$ to $a_{n+1}$ lies in the group
$\Hom_R(\Omega(R/\Ss),M[1]) \cong \Hom_k(\Omega(R/\Ss) \otimes_{\Ss}
k,M[1])$, and since $\Omega(R/\Ss) \otimes_{\Ss} k$ is projective
and $M$ is connective, this group is trivial.
\endproof

\section{T\"oen theorem.}

Now fix an $E_\infty$-algebra $R \in \Dcomm(\Ss)$, and assume given
some $E_1$-algebra $A \in \Dalg(R)$. Then $A$ itself can be
considered not only as a left $A$-module $A \in \D(A)$, but also as
an $R$-module $A \in \D(R)$ and as the diagonal $A$-bimodule $A \in
\D(A^o \otimes_R A)$, where $A^o$ stands for the opposite
$E_1$-algebra.

\begin{defn}
The algebra $A \in \Dalg(R)$ is {\em proper} resp.\ {\em smooth} if
$A$ is compact as an object in $\D(R)$ resp.\ $\D(A^o \otimes_R A)$
is compact.
\end{defn}

Both smoothness and properness are functorial with respect to $R$,
so that sending $R$ to the enhanced category $\Dalg^{sat}(R)$ of
smooth and proper $E_1$-algebras in $\Dalg(R)$ gives an enhanced
functor
\begin{equation}\label{to.eq}
\Dalg^{sat}:\Dcomm(\Ss) \to \eCat.
\end{equation}
The following beautiful theorem has been essentially proved by
B. T\"oen.

\begin{theorem}\label{to.thm}
The functor \eqref{to.eq} commutes with filtered homotopy colimits.
\end{theorem}

Strictly speaking, T\"oen in \cite{to} only considered the
situations when $R$ is a commutative ring; let us recall the
argument to see that it works for spectral algebras with no changes
whatsoever.

\begin{defn}
Assume given an $E_\infty$-algebra $R \in \Dcomm(\Ss)$ and two
$E_1$-algebras $A,B \in \Dalg(R)$, and denote $\D(A,B) = \D(A
\otimes_R B)$. Then an object $M \in \D(A, B)$ is {\em coherent} if
it is compact as a $B$-module.
\end{defn}

T\"oen uses ``pseudoperfect'' instead of ``coherent'' but coherent
is shorter. It is also consistent with earlier terminology: for any
$A \in \Dalg(R)$, we have $\D(A,R) = \D(A \otimes_R R) = \D(A)$, and
this identification identifies coherent objects. For any $A,B \in
\Dalg$, we denote by $\D(A,B)^{coh} \subset \D(A,B)$ the full
enhanced subcategory spanned by coherent objects. We observe that
for any $A$, the diagonal bimodule $A \in \D(A^o,A) = \D(A^o
\otimes_R A)$ is always coherent.

\begin{lemma}\label{pspf.le}
An $E_1$-algebra $A \in \Dalg(R)$ is smooth resp.\ proper if and
only if for any $B \in \Dalg(R)$, $\D(A^o,B)^{coh} \subset
\D(A^o,B)^{pf}$ resp.\ $\D(A^o,B)^{pf} \subset \D(A^o, B)^{coh}$.
\end{lemma}

\proof{} For properness, note that the free right $A$-module $A \in
\D(A^o)$ is compact, so that if $\D(A^o)^{pf} \subset
\D(A^o)^{coh}$, then $A$ is coherent --- that is, compact over
$R$. Conversely, being coherent is closed under retracts and finite
homotopy colimits, so it suffices to check that if $A$ is proper,
then $A^o \otimes_R S \otimes_R B$ is coherent for any compact $S
\in \D(R)$, and this is obvious.

For smoothness, recall that $A \in \D(A^o,A)$ lies in $\D(A^o
\otimes_R A)^{coh}$, so that if $\D(A^o,A)^{coh} \subset
\D(A^o,A)^{pf}$, it is compact. Conversely, note that for any $B$,
any coherent $M \in \D(A^o,B)$ and any compact $N \in \D(A^o,A)$, $N
\otimes_A M \in \D(A^o,B)$ is compact --- indeed, it again suffices
to check this for $N = A^o \otimes_R S \otimes_R A$ for some compact
$S \in \D(R)$, and then $N \otimes_A M \cong A \otimes_R S \otimes_R
M$. But then if $A$ is smooth, any coherent $M \cong A \otimes_A M$
in $\D(A^o, B)$ is therefore compact.
\endproof

\begin{lemma}\label{sat.le}
A smooth and proper $E_1$-algebra $A \in \Dalg(R)$ is compact.
\end{lemma}

\proof{} For any two algebras $A,B \in \Dalg(R)$, the $\Hom$-space
$\Hom(A,B)$ of the enhanced category $\Dalg(R)$ fits into a
functorial homotopy cartesian square
$$
\begin{CD}
\Hom(A,B) @>>> \Iso(\D(A^o,B)^{coh})\\
@VVV @VVV\\
\ppt @>>> \Iso(\D(B)^{pf}),
\end{CD}
$$
where $\Iso$ stands for the enhanced isomorphism groupoid of an
enhanced category, the rightmost arrow is the forgetful functor, and
the bottom arrow is the embedding onto $B \in \D(B)$. If $A$ is
smooth and proper, we can replace coherent objects with compact
ones by Lemma~\ref{pspf.le}, and then recall that $\D^{pf}$ commutes
with filtered homotopy limits by Proposition~\ref{filt.prop}. Since
filtered homotopy colimits commute with finite homotopy limits,
this proves that $\Hom(A,-)$ also commutes with filtered homotopy
colimits.
\endproof

\begin{lemma}\label{hpf.le}
A compact algebra $A \in \Dalg(R)$ is smooth.
\end{lemma}

\proof{} For any bimodule $M \in \D(A^o,A)$, we have the split
square-zero extension $A \oplus M \in \Dalg(R)$, and its splittings
$A \to A \oplus M$ correspond to maps $I_A \to M$ from a
non-commutative version $I_A \in \D(A^o \otimes_R A)$ of the
cotangent module. This bimodule $I_A$ fits into an exact triangle
$$
\begin{CD}
I_A @>>> A^o \otimes_R A @>>> A @>>>
\end{CD}
$$
in the triangulated category $\pi_0(\D(A^o,A))$, thus it is
compact if and only if so is the diagonal bimodule $A$.
\endproof

\proof[Proof of Theorem~\ref{to.thm}.] By Lemma~\ref{sat.le}, we
have a full embeding $\Dalg^{sat}(R) \subset \Dalg(R)^{pf}$ for any
$R \in \Dcomm(\Ss)$, so that in particular, $\Dalg^{sat}(R)$ is
small, and then by Proposition~\ref{filt.prop}, $\Dalg^{pf}$
commutes with filtered homotopy colimits. Therefore for any enhanced
functor $R_\idot:I \to \Dcomm(\Ss)$ with small filtered $I$ and $R
\cong \hocolim_IR_\idot$, the functor
$$
\hocolim_I\Dalg^{sat}(R_\idot) \to \Dalg^{sat}(R)
$$
is fully faithful, and we only need to check that it is essentially
surjective. In other words, we may assume given an algebra $A_i \in
\Dalg(R_i)^{pf}$ such that $A = A_i \otimes_{R_i} R$ is proper, and
we need to show that for some map $i \to i'$, $A_{i'} = A_i
\otimes_{R_i} R_{i'}$ is already proper (while smoothness is
guaranteed by Lemma~\ref{hpf.le}).

Since $A$ is proper and $\D^{pf}$ commutes with filtered homotopy
colimits, we may assume that as an $R$-module, $A \cong M_{i'}
\otimes_{R_{i'}} R$ for some $i' \in I$ and $M_{i'} \in
\D(R_{i'})$. Since $I$ is filtered, we can choose $i'' \in I$ with
maps $i \to i''$, $i' \to i''$, and then replacing $I$ with $i''
\setminus I$, we may assume that $I$ has an initial object $o$, $A_o
\in \Dalg(R_o)$ is compact, and we have an isomorphism of
$R$-modules $A_o \otimes_{R_o} R \cong M_o \otimes_{R_o} R$ for some
$M_o \in \D(R_o)^{pf}$. Denote $A_i = A_o \otimes_{R_o} R_i$, $M_i =
M_o \otimes_{R_o} R_i$, $i \in I$. Since $A$ is a coherent
$A$-module, its $A$-module structure is induced by restriction via
an action map $a:A \to \End_R(A)$. By restriction, $\End_R(A)$ is an
$R_o$-algebra, and the map $a$ is adjoint to a map $a_o:A_o \to
\End_R(A)$ in $\Dalg(R_o)$. But
$$
\End_R(A) \cong \End_{R_o}(M_o) \otimes_{R_o} R \cong
\hocolim_I\End_{R_i}(M_i),
$$
and since $A_o \in \Dalg(R_o)$ is compact, the map $a_o$ factors
through some map $A_o \to \End_{R_i}(M_i)$, $i \in I$ adjoint to a
map $a_i:A_i \to \End_{R_i}(M_i)$ in $\Dalg(R_i)$. Then by
restriction, $M_i$ becomes a coherent $A_i$-module, and since $A_i$
is compact, $M_i \in \D(A_i)$ is compact by Lemma~\ref{hpf.le} and
Lemma~\ref{pspf.le}. Thus we have two compact $A_i$-modules, $M_i$
and $A_i$ itself, and an isomorphism $A_i \otimes_{R_i} R \cong M_i
\otimes_{R_i} R$ in $\D(A)$. Since $\D^{pf}$ commutes with filtered
homotopy colimits by Proposition~\ref{filt.prop}, this isomorphism
must be induced by an isomorphism $A_{i'} \cong M_{i'}$ in
$\D(A_{i'})$ for some $i' \in I$. But $M_{i'} \in \D(A_{i'})$ is not
only compact but also coherent, so that $A_{i'}$ must be proper.
\endproof

\section{Tate diagonal.}

Recall that for any $R \in \Dcomm(\Ss)$ and any set $S$, we denote
by $R[S]$ the direct sum of copies of $R$ numbered by elements $s
\in S$. More generally, for a topological space $X$, we denote by
$R[X] \in \D(\Ss)$ the $R$-homology spectrum of $X$. If $X=G$ is a
compact Lie group, then $R[G]$ is a $E_1$-algebra in $\Dalg(R)$ with
respect to the Pontryagin product, and the projection $G \to \ppt$
induces the augmentation $E_1$-map $R[G] \to R$. Restricting with
respect to the augmentation gives a tautological embedding $a:\D(R)
\to \D(R[G])$ that has adjoints on the left and on the right, $M
\mapsto M_{hG}$ resp.\ $M \mapsto M^{hG}$, known as the {\em
  homotopy quotient} and the {\em homotopy fixed points}
functors. If $R$ is discrete, thus simply a ring, and the group $G$
is finite, then $\D(R[G])$ is the derived category of $R$-linear
representations of the group $G$, homotopy quotient is group
homology, and homotopy fixed points is group cohomology. In the
general situation, the diagonal embedding $G \to G \times G$ turns
$\D(R[G])$ into a symmetric monoidal enhanced category, the
tautological embedding $a$ is symmetric monoidal, and the homotopy
fixed points functor is lax symmetric monoidal by adjunction. Thus
in particular, $R^{hG}$ is naturally an $E_\infty$-algebra in
$\D(R)$, and the homotopy fixed points functor can be refined to a
functor
\begin{equation}\label{hG.rf}
\D(R[G]) \to \D(R^{hG}), \qquad M \mapsto M^{hG}.
\end{equation}
Since $G$ is assumed to be compact, the algebra $R[G]$ is proper, so
that we have a full embedding
\begin{equation}\label{pf.coh}
\D(R[G])^{pf} \subset \D(R[G])^{coh}
\end{equation}
and the induced embedding
\begin{equation}\label{ind.pf.coh}
\Ind(\D(R[G])^{pf}) = \D(R[G]) \subset \D(R[G])^{coh}.
\end{equation}
However, $R[G]$ is usually not smooth, so that the embeddings
\eqref{pf.coh}, \eqref{ind.pf.coh} are not equivalences. We then
have a non-trivial enhanced Verdier quotient
$$
\D(R[G])^{sing} = \D(R[G])^{coh}/\D(R[G])^{pf}.
$$
The subcategory $\D(R[G])^{coh} \subset \D(R[G])$ is symmetric
monoidal, and the subcategory $\D(R[G])^{pf} \subset \D(R[G])^{coh}$
is a symmetric monoidal ideal, so that $\D^{sing}(R[G])$ is also a
symmetric monoidal enhanced category in a natural way. On the level
of $\Ind$-completions, \eqref{ind.pf.coh} induces a semiorthogonal
decomposition
\begin{equation}\label{deco}
\Ind(\D(R[G])^{coh}) = \langle \Ind(\D(R[G])^{sing}),\D(R[G])
\rangle.
\end{equation}
The stable enhanced categories $\Ind(\D(R[G])^{coh})$,
$\Ind(\D(R[G])^{sing})$ are symmetric monoidal, and so is the
projection
\begin{equation}\label{l.eq}
l:\Ind(\D(R[G])^{coh}) \to \Ind(\D(R[G])^{sing})
\end{equation}
onto the first factor of the decomposition \eqref{deco}. The
augmentation functor $\D(R) \to \D(R[G])$ composed with
\eqref{ind.pf.coh} and the projection $l$ provides a symmetric
monoidal functor $\D(R) \to \Ind(\D(R[G])^{sing})$ that has a
right-adjoint {\em Tate fixed points} functor
\begin{equation}\label{tG}
\Ind(\D(R[G])^{sing}) \to \D(R), \qquad M \mapsto M^{tG}.
\end{equation}
By abuse of notation, we will denote $M^{tG} = l(M)^{tG}$ for any
$M$ in the category $\Ind(\D(R[G])^{coh}$ (and in particular, for
any coherent $M \in \D(R[G])$). The functor \eqref{tG} is lax
symmetric monoidal by adjunction, so that $R^{tG}$ is an
$E_\infty$-algebra in $\D(R)$, and then as in \eqref{hG.rf},
\eqref{tG} can be canonically refined to an enhanced functor
\begin{equation}\label{tG.rf}
\Ind(\D(R[G])^{sing}) \to \D(R^{tG}), \qquad M \mapsto M^{tG}.
\end{equation}
For any $M \in \D(R[G])^{coh}$, the decomposition \eqref{deco}
induces an exact triangle
\begin{equation}\label{deco.tri}
\begin{CD}
M_{hG}[d] @>{t}>> M^{hG} @>>> M^{tG} @>>>,
\end{CD}
\end{equation}
where $d = \dim G$ is the dimension of $G$, and $t$ is a natural
trace map induced by the Poincar\'e duality on $G$ (if the group $G$ is
finite, $t$ is just the averaging over the group).

Sometimes Tate fixed points can be computed by localizing the usual
homotopy fixed points with respect to certain elements in the
homotopy groups of $R^{hG}$. The basic example is $G = S^1$, the
unit circle. If (and only if) $R \in \Dcomm(\Ss)$ is orientable as a
multiplicative generalized cohomology theory --- for example, if $R$
is discrete --- we have $\pi_\idot(R^{hS^1}) \cong \pi_\idot(R)[u]$,
where $u$ is a single generator of cohomological degree $2$. In this
case, $\pi_\idot(R^{tG}) = \pi_\idot(R)[u,u^{-1}]$, and for any $M
\in \D(R[S^1])^{coh}$, we have
\begin{equation}\label{S.1.u}
M^{tS^1} \cong M^{hS^1} \otimes_{R^{hS^1}} R^{tS^1} \cong \hocolim_n
M^{tS^1}[2n],
\end{equation}
where the colimit is taken with respect to the action $u:M^{hS^1}
\to M^{hS^1}[2]$ of the generator $u \in \pi_{-2}(R^{hS^1})$.

Another example is when $G=C_p \subset S^1$ is the cyclic group of
some prime order $p \geq 3$, and $R$ is a ring annihilated by
$p$. In this case, $\pi_\idot(R^{hC_p}) \cong R\langle \eps,u
\rangle$, where $u$ has cohomological degree $2$, $\eps$ has
cohomological degree $1$, and they commute. Tate fixed points
$R^{tC_p}$ are again obtained by inverting $u$, and for any coherent
$M \in \D(R[C_p])$, we again have
\begin{equation}\label{c.p.u}
M^{tC_p} \cong \hocolim_nM^{hC_p}[2n],
\end{equation}
with colimit takes with respect to the action on $u$.

If $R$ is not orientable, $\pi_\idot(R^{tS^1})$ is still the
abutment of an Atiyah-Hizberuch spectral sequence whose first page
is $\pi_\idot(R)[u]$, but the spectral sequence does not degenerate,
and the periodicity element $u$ does not survive to the last
page. We do not know any general method to compute $R^{tS^1}$. The
situation for the cyclic group is similar; however, there is the
following striking result.

\begin{lemma}\label{tate.le}
Let $R=\Ss_q$, $q=p^n$ be the $n$-fold \'etale covering of the
$p$-completion $\Ss_p$ of the sphere, for some $n \geq 1$ and some
prime $p$. For any $M \in \D(R)$, consider $M^{\otimes_R p}$ as an
object in $\D(R[C_p])$ via the longest cycle permutation
action. Then there is a map
\begin{equation}\label{td}
M \to (M^{\otimes_R p})^{tC_p},
\end{equation}
functorial in $M$, and this map is an isomorphism if $M$ is
compact.\endproof
\end{lemma}

This is a version of the Segal Conjecture, see \cite[III.1]{NS} and
references therein. Nikolaus and Scholze call \eqref{td} the {\em
  Tate diagonal map}. The essential part of the proof is the case $M
= R$ (when $M^{\otimes_R p}$ is again $R$ with the trivial
$C_p$-action).

Our proof of Hodge-to-de Rham Degeneration relies on one immediate
corollary of Lemma~\ref{tate.le}. Observe that for any map $R_0 \to
R_1$ in $\Dcomm(\Ss)$, the augmentation embedding commutes with the
tensor product functor $- \otimes_{R_0} R_1$, so that by adjunction, we
obtain a functorial map
\begin{equation}\label{bc.t}
M^{tG} \otimes_{R_0^{tG}} R_1^{tG} \to (M \otimes_{R_0} R_1)^{tG}
\end{equation}
for any coherent $M$ in $\D(R_0[G])$, where $M^{tG}$ is considered
as an $R_0^{tG}$-module via the refinement \eqref{tG.rf}.

\begin{corr}\label{tate.corr}
Let $R$ be as in Lemma~\ref{tate.le}, and let $k = \F_q$, $q = p^n$
be the degree-$n$ Galois extension of the prime field $\F_p$, with
the natural map $R \to k$. Then for any compact $M \in \D(R)$ with
$M_k = M \otimes_R k$, the map
\begin{equation}\label{tate.iso}
M_k \otimes_k k^{tC_p} \to (M^{\otimes_k p})^{tC_p}
\end{equation}
obtained by composing \eqref{td} and \eqref{bc.t} is an isomorphism.
\end{corr}

\proof{} Both sides are functorial in $M$, and the functors are
stable enhanced functors, thus commute with finite homotopy colimits
and with retracts. Therefore it suffices to consider the case $M=R$
where the statement immediately follows from Lemma~\ref{tate.le}.
\endproof

\section{Hodge-to-de Rham degeneration.}

\subsection{Cyclic homology.}

For any $E_\infty$-algebra $R \in \Dcomm(\Ss)$ and any $E_1$-algebra
$A \in \Dalg(R)$ over $R$, the {\em Hochschild Homology} of $A$ over
$R$ is defined as the $R$-module $HH(A/R) = A^o \otimes_{A^o
  \otimes_R A} A$. To describe it more explicitly, one uses the bar
construction to replace $A$ with a termwise-free simplicial
$A$-bi\-module; this provides a canonical enhanced functor
$(A/R)^\Delta_\hash:\Delta^o \to \D(R)$ and an identification
$$
HH(A/R) \cong \hocolim_{\Delta^o} (A/R)^\Delta_\hash.
$$
It is well-known that $HH(A/R)$ can be promoted to an object in
$\D(R[S^1])$. To construct the $S^1$-action, one observes that
$(A/R)^\Delta_\hash$ extends to A. Connes' cyclic category $\Lambda$
of \cite{connes}: we have an embedding $j:\Delta^o \to \Lambda$ and
an enhanced functor $(A/R)_\hash:\Lambda \to \D(R)$ such that
$j^*(A/R)_\hash \cong (A/R)_\hash^\Delta$. For any enhanced functor
$E:\Lambda \to \D(R)$, one defines
$$
HH(E) = \hocolim_{\Delta^o}j^*E, \qquad HC(E) = \hocolim_\Lambda E,
$$
and one proves that $HH$ extends to a functor $\wt{HH}:\D(R)^\Lambda
\to \D(R[S^1])$ (the cleanest construction of this extension is
given in \cite{dr}). In fact, one can say more: the classifying
space $|\Lambda|$ of the nerve of the category $\Lambda$ is
canonically identified with the classifying space $BS^1$ of the
circle, and $\D(R[S^1])$ is naturally identified with the full
subcategory in $\D(R)^\Lambda$ spanned by locally constant enhanced
functors. The functor $\wt{HH}$ is left-adjoint to the full embedding
$\D(R[S^1]) \subset \D(R)^\Lambda$. This implies that $HC(E) \cong
HH(E)_{hS^1}$, and this is known as cyclic homology. Periodic cyclic
homology $HP(E)$ is then defined as
$$
HP(E) = HH(E)^{tS^1},
$$
and one shortens $HP((A/R)_\hash)$, $HC((A/R)_\hash)$ to $HP(A/R)$
resp.\ $HC(A/R)$. If $R$ is discrete, thus oriented, then $R^{hS^1}
\cong R[u]$, $R^{tS^1} \cong R[u,u^{-1}]$, and for any $A \in
\Dalg(R)$ we have spectral sequences
\begin{equation}\label{hdr}
HH(A/R)[u^{-1}] \Rightarrow HC(A/R), \qquad HH(A/R)((u)) \Rightarrow
HP(A/R).
\end{equation}
These are known as the {\em Hodge-to-de Rham spectral sequences}.

For any integer $n \geq 1$, we have the cyclic subgroup $C_n \subset
S^1$, and its action on $HH$ can be seen directly in terms of the
category $\Lambda$. To do this, one defines a category $\Lambda_n$
equipped with an {\em edgewise subdivision functor} $i_n:\Lambda_n
\to \Lambda$ and a projection $\pi_n:\Lambda_n \to \Lambda$. The
projection $\pi_n$ is a bifibration in groupoids whose fiber $\ppt_n
= \ppt/C_n$ is the connected groupoid with a single object with
automorphism group $C_n$. On the level of classifying spaces,
$|i_n|:|\Lambda_n| \to |\Lambda|$ is a homotopy equivalence, and the
fibration $|\pi_n|:|\Lambda_n| \cong |\Lambda| \to |\Lambda|$ is
obtained by delooping once the short exact sequence
\begin{equation}\label{c.p}
\begin{CD}
1 @>>> C_n @>>> S^1 @>>> S^1 @>>> 1
\end{CD}
\end{equation}
of abelian compact Lie groups. The embedding $j:\Delta^o \to
\Lambda$ fits into a commutative diagram
$$
\begin{CD}
\Delta^o @<{\overline{\pi}_n}<< \Delta^o \times \ppt_n
@>>> \Delta^o\\
@V{j}VV @VV{j_n}V @VV{j}V\\
\Lambda @<{\pi_n}<< \Lambda_n @>{i_n}>> \Lambda,
\end{CD}
$$
where the square on the left is cartesian, and
$\overline{\pi}_n:\Delta^o \times \ppt_n \to \Delta^o$ is the
projection onto the first factors. The classical Edgewise
Subdivision Lemma \cite{edge} shows that for any $E \in
\D(R)^\Lambda$, the natural map
\begin{equation}\label{edge.eq}
\hocolim_{\Delta^o}j_n^*i_n^*E \to \hocolim_{\Delta^o}j^E = HH(E)
\end{equation}
is an isomorphism, and its source lies naturally in $\D(R)^{\ppt_n}
\cong \D(R[C_n])$.

This construction is especially useful if $n=p$ is an odd prime, and
$R$ is a ring annihilated by $p$. Namely, for any $R$ and $M \in
\D(R[S^1])$, the exact sequence \eqref{c.p} provides an
identification
\begin{equation}\label{h.c.p}
(M^{hC_p})^{hS^1} \cong M^{hS^1}.
\end{equation}
If $R$ is a ring annihilated by $p$, then the left-hand side carries
two periodicity endomorphisms of degree $2$: $u$ coming from $C_p$,
and $u'$ coming from $S^1 = S^1/C_p$. Hochschild-Serre spectral
sequence for \eqref{c.p} shows that it is the first endomorphism $u$
that is compatible with the periodicity endomorphism $u$ in the
right-hand side, so that \eqref{h.c.p} coupled with \eqref{S.1.u}
and \eqref{c.p.u} provides a map
\begin{equation}\label{t.c.p}
M^{tS^1} \to (M^{tC_p})^{hS^1}.
\end{equation}
Moreover, the same Hochschild-Serre spectral sequence shows that
$u'$ actually vanishes, so that $(M^{tC_p})^{tS^1}=0$, and we have
$(M^{tC_p})^{hS^1} \cong (M^{tC_p})_{hS^1}[1]$ by
\eqref{deco.tri}. Since homotopy quotients commute with homotopy
colimits, we conclude that \eqref{t.c.p} is an isomorphism. This
allows one to reduce questions about $M^{tS^1}$ to questions about
$M^{tC_p}$.

\subsection{Degeneration theorem.}

We can now state and prove the Hodge-to-de Rham Degeneration
Theorem. First, assume given a ring $k$ annihilated by an odd prime
$p$, and an algebra $A \in \Dalg(k)$. Consider the corresponding
enhanced functor $(A/k)_\hash:\Lambda \to \D(k)$ and its edgewise
subdivision $j_n^*i_n^*(A/k)_\hash$of \eqref{edge.eq}. We then have
natural map
\begin{equation}\label{hoco.d}
\hocolim_{\Delta^o}(j_n^*i_n^*(A/k)_\hash)^{tC_p} \to
\left(\hocolim_{\Delta^o}(j_n^*i_n^*(A/k)_\hash)\right)^{tC_p}
\end{equation}
in $\D(k[S^1])$, and its target is identified with $HH(A/k)^{tC_p}$
by \eqref{edge.eq}.

\begin{lemma}\label{sm.le}
Assume that the algebra $A$ is smooth. Then the map \eqref{hoco.d}
is an isomorphism.
\end{lemma}

\proof{} Since $A$ is smooth, the diagonal bimodule $A$ is compact,
and then it is a retract of some piece of the stupid filtration of
its bar resolution. Therefore for some $n \geq 1$, the homotopy
colimits $\hocolim_{\Delta^o}$ in \eqref{hoco.d} are retracts of
homotopy colimits $\hocolim_{\Delta^o_{\leq n}}$ over the full
subcategory $\Delta^o_{\leq n} \subset \Delta^o$ spanned by finite
totally ordered sets with at most $n$ elements. But the category
$\Delta_{\leq n}^o$ is finite, and the Tate fixed ponts functor
$(-)^{tC_p}$, being stable, commutes with finite homotopy colimits.
\endproof

\begin{remark}
If $A$ is not smooth, \eqref{hoco.d} is not an isomorphism, but its
source still has an invariant meaning --- in fact,
$\hocolim_{\Lambda}(i_n^*(A/k)_\hash)^{tC_p}$ is the so-called {\em
  co-periodic cyclic homology} $\overline{HP}(A/k)$, a new
localizing invariant of DG-algebras introduced and studied in
\cite{coper}. Mathew in \cite{M} has no counterpart of
Lemma~\ref{sm.le}, and co-periodic cyclic homology does not appear
explicitly. It seems that the real reason for this is that he uses
Topological Hochschild Homology $THH(A)$, and one can show that for
a DG algebra $A$ over a finite field $k$, $THH(A)$ becomes
isomorphic to $\overline{HP}(A/k)$ after one inverts the B\"okstedt
periodicity generator. We will return to this elsewhere.
\end{remark}

Next, let $k = \F_q$, $q = p^n$ be a finite field of odd
characteristic $p$, and let $R = \Ss_q$ be as in
Corollary~\ref{tate.corr}.

\begin{lemma}\label{hp.le}
Assume given a smooth and proper algebra $A \in \D(R)$, with $A_k =
A \otimes_R k$. Then there exists an isomorphism
\begin{equation}\label{hp.i}
HP(A/k) \cong HH(A/k) \otimes_k k[u,u^{-1}].
\end{equation}
\end{lemma}

\proof{} Consider the enhanced functor $(A/R)_\hash:\Lambda \to
\D(R)$, its edgewise subdivision $i_p^*(A/R)_\hash$, and its
restriction $j_p^*i_p^*(A/R)_\hash$ to $\Delta^o \times \ppt_p
\subset \Lambda_p$. We have a natural identification
$j_p^*i_p^*(A/R)_\hash \cong
(\overline{\pi}_p^*j^*(A/R)_\hash)^{\otimes_R p}$, and the Tate
diagonal map \eqref{td} then induces a map $j^*(A/R)_\hash \to
(j_p^*i_p^*(A/R)_\hash)^{tC_p}$. Moreover, it has been shown in
\cite[III.2]{NS} that this map extends to a map
\begin{equation}\label{tl}
(A/R)_\hash \to \pi_{p*}^ti_p^*(A/R)_\hash,
\end{equation}
where $\pi_{n*}^t:\D(R)^{\Lambda_p} \to \D(R)^\Lambda$ is the
relative version of the Tate fixed points functor for the
bifibration $\pi_p:\Lambda_p \to \Lambda$. Then as in
Corollary~\ref{tate.corr}, the map \eqref{tl} induces a map
\begin{equation}\label{tl.k}
(A_k/k)_\hash \otimes_k \pi_{p*}^t \to \pi_{p*}i^*_p(A_k/k)_\hash,
\end{equation}
and since $A_k$ is proper, this map is an isomorphism. But $A_k$ is
also smooth, and then by Lemma~\ref{sm.le}, the isomorphisms
\eqref{t.c.p} and \eqref{hoco.d} provide an isomorphism
\eqref{hp.i}.
\endproof

\begin{remark}
We note that one does not need the full force of Lemma~\ref{tate.le}
to obtain Corollary~\ref{tate.corr} and Lemma~\ref{hp.le}. In
effect, for any complex $M_\idot$ of $k$-vector spaces, one can
equip $(M_\idot^{\otimes_k p})^{tC_p}$ with a natural
$C_p$-equivariant $\Z$-indexed increasing filtration $\beta_\idot$
whose associated graded quotients $\gr^{\beta}_n$ are the shifts
$M_\idot[n]$, and the quotient map $\beta_0(M_\idot^{\otimes_k
  p})^{tC_p} \to M_\idot$ admits a canonical $\Ss$-linear
splitting. If $M_\idot$ is of the form $M_\idot = M \otimes_{\Ss} k$
for a spectrum $M$, the splitting can be made $k$-linear, and this
provides isomorphisms \eqref{tate.iso} and \eqref{tl.k}. This is the
approach taken explicitly in \cite{Ka1} and implicitly in \cite{Ka2}
(where the spectrum $M$ is not mentioned by name, and the only thing
used are obstructions to its existence). The construction using
Lemma~\ref{tate.le} is obviously much more direct and conceptually
clear, but this comes at a price: we have to use the proof of the
Segal Conjecture as a black box. It would be interesting to see if
the technology of \cite{coper} and \cite{Ka2} can clarify the
contents of the black box.
\end{remark}

\begin{theorem}
Assume given a smooth and proper algebra $A \in \Dalg(K)$ over a
field $K$ of characteristic $0$. Then the Hodge-to-de Rham spectral
sequence for $HP(A/K)$ degenerates.
\end{theorem}

\proof{} By Proposition~\ref{sm.prop} and Theorem~\ref{to.thm}, one
can choose a formally smooth $E_\infty$-algebra $R \in \Dcomm(\Ss)$
equipped with a map $a:R \to K$, and a smooth and proper algebra
$A_R \in \Dalg(R)$ such that $A_R \otimes_R K \cong A$. Localizing
$R$ if necessary, we may assume that it lies in
$\Dcomm(\Ss(2^{-1}))$. The map $a$ factors through the finitely
generated ring $\bR = \pi_0(R)$, and if we let $A_{\bR} = A_R
\otimes_R \bR$, then it suffices to prove that the Hodge-to-de Rham
spectral sequence for $HP(A_{\bR}/\bR)$ degenerates. Since $A_R$ is
smooth and proper, $A_{\bR}$ is also smooth and proper, so that
Hochschild homology groups $HH_\idot(A_{\bR}/\bR)$ are finitely
generated $\bR$-modules. Then by Nakayama Lemma, to prove that all
the differentials in the spectral sequence vanish, it suffices to
prove that for any residue field $k$ of the ring $\bR$, with $A_k =
A_{\bR} \otimes_{\bR} k$, the Hodge-to-de Rham spectral sequence
$HP(A_k/k)$ degenerates. But this is a spectral sequence of
finite-dimensional $k$-vector spaces, $k$ is a finite field of odd
characteristic, and by Lemma~\ref{sm.le} and Lemma~\ref{fm.le}, its
first and last page have the same dimensions,
\endproof

{\small\noindent {\em Affiliations\/}: A.K., K.M.: National Research
  University Higher School of Economics (NRU HSE), D.K.: Steklov
  Mathematics Institute and NRU HSE.

\noindent
{\em E-mail addresses\/}: {\tt kaledin@mi.ras.ru} (D.K.), {\tt
  kon\_an\_litsey@list.ru} (A.K.), {\tt kirill.salmi94@gmail.com}
(K.M.)}

\end{document}